\documentclass[12pt]{amsart}
\usepackage{amssymb,latexsym,epsfig}

\input{epsf}
\topskip  \numberwithin{equation}{section}
\newtheorem{theorem}{Theorem}[section]
\newtheorem{proposition}[theorem]{Proposition}
\newtheorem{corollary}[theorem]{Corollary}

\newtheorem{lemma}[theorem]{Lemma}

\begin{document}
\def\sof{\hfill\rule{2mm}{2mm}}
\def\ls{\leq}
\def\gs{\geq}
\def\PP{\mathcal P}
\def\qq{{\bold q}}
\def\txx{{\frac1{2\sqrt{x}}}}

\title{Enumeration of $(k,2)$-noncrossing partitions}
\author{}
\maketitle

\pagenumbering{arabic} \pagestyle{headings}

\begin{center}
Toufik Mansour$^{1}$ and Simone Severini$^{2}$

$^{1}$Department of Mathematics, University of Haifa, 31905 Haifa, Israel.

$^{1}$Center for Combinatorics, LPMC, Nankai University, 300071 Tianjin, P.
R. China

$^{2}$Department of Mathematics and Department of Computer Science,
University of York, YO10 5DD York, United Kingdom

\texttt{toufik@math.haifa.ac.il, ss54@york.ac.uk}
\end{center}


\section*{Abstract}

A set partition is said to be $(k,d)$-noncrossing if it avoids the
pattern $12\cdots k12\cdots d$. We find an explicit formula for
the ordinary generating function of the number of
$(k,d)$-noncrossing partitions of $\{1,2,\ldots,n\}$ when $d=1,2$.
\medskip

\noindent {Keywords}: partitions, forbidden subsequences, kernel method.

\noindent \textsc{2000 Mathematics Subject Classification}: 05A05, 05A15

\section{Introduction}

A \emph{partition} $\Pi $ of the set $[n]=\{1,2,\ldots ,n\}$ is a collection
$B_{1},B_{2},\ldots ,B_{d}$ of nonempty disjoint subsets of $[n]$. The
elements of a partition are called \emph{blocks}. We assume that $%
B_{1},B_{2},\ldots ,B_{d}$ are listed in the increasing order of their
minimum elements, that is $\min B_{1}<\min B_{2}<\cdots <\min B_{d}$. The
set of all partitions of $[n]$ with $d$ blocks is denoted by $P(n,d)$. The
cardinality of $P(n,d)$ is the well-known Stirling number of the second kind
\cite{S1}, which is usually denoted by $S(n,k)$. Any partition $\Pi $ can be
written in the \emph{canonical sequential form} $\pi _{1}\pi _{2}\cdots \pi
_{n}$, where $i\in B_{\pi _{i}}$ (see, \emph{e.g.} \cite{K1}). From now on,
we identify each partition with its canonical sequential form. For example,
if $\Pi =\{1,4\},\{2,5,7\},\{3\},\{6\}$ is a partition of $[7]$, then its
canonical sequential form is $\pi =1231242$ and in such a case we write $\Pi
=\pi $.

The \emph{reduced form} of a word $\pi $ on the alphabet $%
\{a_{1},a_{2},\ldots ,a_{d}\}$, where $a_{1}<a_{2}<\cdots <a_{d}$, is a word
$\pi ^{\prime }$ obtained by renaming the letters of $\pi $. Specifically,
the letter $a_{i}$ is renamed $i$ for all $i=1,2,\ldots ,d$. For example,
the reduced form of the word $135351$ is $123231$. We say that a partition $%
\Pi $ \emph{avoids} $\tau =\tau _{1}\cdots \tau _{k}$ if its canonical
sequential form $\pi =\pi _{1}\pi _{2}\cdots \pi _{n}$ does not contain a
subsequence $\pi ^{\prime }=\pi _{i_{1}}\cdots \pi _{i_{k}}$, such that the
reduced form of $\pi ^{\prime }$ equals the reduced form of $\tau $. For
example, the partition $\pi =1231242$ avoids $12321$.

In \cite{K1,K2}, Klazar proved that the number of \emph{noncrossing
partitions} of $[n]$, that is partitions avoiding $1212$, and the number of
\emph{nonnesting partitions} of $[n],$that is partitions avoiding $1221$,
are given by the $n$-th Catalan number $\frac{1}{n+1}\binom{2n}{n}$. These
results of Klazar have been extended in several recent directions. For
instance, Chen\emph{\ et al.} \cite{C1, C2} studied $m$-regular noncrossing,
$k$-noncrossing and $k$-nonnesting partitions. Recall that a partition $\Pi $
is called $m$\emph{-regular} if $|x_{1}-x_{2}|\geq m$, for any two distinct
elements $x_{1},x_{2}$ in the same block of $\Pi $. A partition $\pi $ is
called $k$\emph{-noncrossing} if it avoids $12\cdots k12\cdots k$ and $k$%
\emph{-nonnesting} if it avoids $12\cdots kk\cdots 21$. For
further works on this subject, the reader is referred to Sagan
\cite{SA} and the references therein.

In this paper, we generalize the concept of $k$-noncrossing
partitions to $(k,d)$-noncrossing partitions. A partition is said
to be $(k,d)$\emph{-noncrossing} if it avoids $12\cdots k12\cdots
d$. Let $\mathcal{N}_{k,d}(n)$ be the set of all
$(k,d)$-noncrossing partitions of $[n]$. Note that
$\mathcal{N}_{k,k}(n)$ is the set of $k$-noncrossing partitions of
$[n]$ (see \cite{C2}).

For $d=0$, it is easy to see from the definitions that the number
of $(k,0)$-noncrossing partitions of $[n]$ is the same as the
number of partitions of $[n]$ with at most $k-1$ blocks. Thus,
\begin{equation*}
\#\mathcal{N}_{k,0}(n)=\sum_{i=0}^{k-1}S(n,i),
\end{equation*}
where $S(n,i)$ is the Stirling number of the second kind. In this
paper we give a complete enumeration for two cases of
$(n,d)$-noncrossing partitions, in which $d$ is either $1$ or $2$.

\begin{center}
\begin{table}[th]
\begin{tabular}{c|ccccccccccccc}
$k\backslash n$ & 0 & 1 & 2 & 3 & 4 & 5 & 6 & 7 & 8 & 9 & 10 & 11 & 12 \\
\hline
2 & 1 & 1 & 2 & 5 & 14 & 42 & 132 & 429 & 1430 & 4862 & 16796 & 58786 &
208012 \\
3 & 1 & 1 & 2 & 5 & 15 & 51 & 188 & 731 & 2950 & 12235 & 51822 & 223191 &
974427 \\
4 & 1 & 1 & 2 & 5 & 15 & 52 & 202 & 856 & 3868 & 18313 & 89711 & 450825 &
2310453 \\
5 & 1 & 1 & 2 & 5 & 15 & 52 & 203 & 876 & 4112 & 20679 & 109853 & 608996 &
3488806 \\
6 & 1 & 1 & 2 & 5 & 15 & 52 & 203 & 877 & 4139 & 21111 & 115219 & 666388 &
4045991%
\end{tabular}%
\label{Table}\vspace{5pt}
\caption{Number the $(k,2)$-noncrossing partitions of $[n]$ for $k=2,3,4,5,6$
and $n=0,1,\ldots ,12$.}
\end{table}
\end{center}

\vspace{-35pt}Table~\ref{Table} presents the number of
$(k,2)$-noncrossing partitions in
$\mathcal{P}_{k}(n)=\mathcal{N}_{k,2}(n)$, where $k=2,3,4,5,6$. We
will show that the ordinary generating function
\begin{equation*}
\sum_{n\geq 0}\#\mathcal{P}_{k}(n)x^{n}
\end{equation*}%
for the number of partitions in $\mathcal{P}_{k}(n)$ is rational
in $x$ and $\sqrt{(1-kx)^{2}-4x^{2}}$. Namely, we prove the
following result.

\begin{theorem}
\label{mth} Let $k\geq 2$ and let
\begin{equation*}
y_{k}=\frac{1-(k-2)x-\sqrt{(1-kx)^{2}-4x^{2}}}{2x(1-(k-2)x)}.
\end{equation*}%
Then the ordinary generating function for the number of $(k,2)$-noncrossing
partitions of $[n]$ is given by
\begin{equation*}
\sum_{n\geq 0}\#\mathcal{P}_{k}(n)x^{n}=\frac{\frac{x^{k-1}y_{k}}{1-xy_{k}}%
+\sum_{j=0}^{k-2}\sum_{i=0}^{j}(-1)^{i+j}x^{i}\beta _{i,j}}{%
1-\sum_{j=0}^{k-2}\sum_{i=0}^{j}(-1)^{i+j}ix\beta _{i,j}},
\end{equation*}%
where $\beta _{j,j}=1$ and
\begin{equation*}
\beta _{i,j}=jx\prod_{s=i+1}^{j-1}(sx-1)
\end{equation*}%
for $i=0,1,\ldots ,j-1$.
\end{theorem}

The proof of Theorem~\ref{mth}, see Section~2, is based on looking at $%
\mathcal{P}_{k}(n)$ as a disjoint union of subsets $\mathcal{P}_{k,\ell }(n)$%
, depending on some parameter $\ell $. Then, we obtain linear recurrence
relations with two indices, $n$ and $\ell $, for the number of partitions in
these subsets. In order to solve the recurrence relations, we make use of
the \emph{kernel method technique} (see, \emph{e.g.}, \cite{B02}). The
subsets $\mathcal{P}_{k,\ell }(n)$ can be defined as the collection of all $%
\ell $-increasing partitions in $\mathcal{P}_{k}(n)$. Recall that $\pi =\pi
_{1}\cdots \pi _{n}$ is an \emph{$\ell $-increasing partition} of $[n]$ if $%
\pi _{i}=i$ for all $i=1,2,\ldots ,\ell $ and
$\tau_{\ell+1}\neq\ell+1$. Let us denote by $\mathcal{Q}_{\ell
}(n)$ the set of all $m$-increasing partitions of $[n]$, with
$m\leq \ell $. Directly from the proof of Theorem~\ref{mth}, we
can obtain a formula for the number of partitions of $[n]$ in
$\mathcal{Q}_{\ell }(n)$.

\begin{corollary}
\label{se}Let $\ell \geq 0$. The ordinary generating function for the number
of partitions in $\mathcal{Q}_{\ell }(n)$ is
\begin{equation*}
I_{\ell }(x)=\sum_{n\geq 0}\#\mathcal{Q}_{\ell }(n)x^{n}=\frac{
1+\sum_{j=1}^{\ell }\sum_{i=0}^{j}(-1)^{i+j}x^{i}\beta _{i,j}}{
1-\sum_{j=1}^{\ell }\sum_{i=0}^{j}\left( -1\right) ^{i+j}ix\beta
_{i,j}},
\end{equation*}
where $\beta _{i,j}$ is defined in Theorem \ref{mth}.
\end{corollary}

\bigskip Theorem~\ref{mth} gives two particular results, namely $k=2$ and $k=3$. For $k=2$,
\begin{equation*}
\sum_{n\geq
0}\#\mathcal{P}_{2}(n)x^{n}=1+\frac{xy_{2}}{1-xy_{2}}=\frac{1}{
1-xy_{2}}=y_{2},
\end{equation*}
where $y_{2}=\frac{1-\sqrt{1-4x}}{2x}$. Thus, the number of
partitions in $ \mathcal{P}_{2}(n)$ is given by the $n$-th Catalan
number. For $k=3$,
\begin{equation*}
\sum_{n\geq
0}\#\mathcal{P}_{3}(n)x^{n}=\frac{1+\frac{x^{2}y_{3}}{1-xy_{3}}}{
1-x}=\frac{3-3x-\sqrt{1-6x+5x^{2}}}{2\left( 1-x\right) },
\end{equation*}
where $y_{3}=\frac{1-x-\sqrt{1-6x+5x^{2}}}{2x\left( 1-x\right) }$.
It follows that the number of partitions in $\mathcal{P}_{3}(n)$
is given by the $n$-th binomial transform of the Catalan number
$\sum_{i=0}^{n}(-1)^{i}3^{n-i}\binom{n}{i}\binom{i}{\lfloor
i/2\rfloor }$ (see \cite[Sequence A007317]{JS}). For $k=2$, there
is a combinatorial proof that the number of partitions in
$\mathcal{P}_{2}(n)$ is given by the $n$-th Catalan number. For
$k=3$, the formula above counts the number of Schr\"{o}der paths
with no peaks at even level of length $n$ (see \cite[Sequence
A007317]{JS}). It would be interesting to find a bijective proof
of this result.

Another bonus from the proof of Theorem \ref{mth} is the ordinary generating
function for the number of $(k,1)$-noncrossing partition of $[n]$.
Specifically, we prove the following result.

\begin{theorem}
\label{mth22} Let $k\geq 2$. Then the ordinary generating function
for the number of $(k,1)$-noncrossing partitions of $[n]$ is given
by
\begin{equation*}
\sum_{n\geq 0}\#\mathcal{N}_{k,1}(n)x^{n}=\frac{1-x+(1-x)\sum_{j=1}^{k-2}%
\sum_{i=0}^{j}(-1)^{i+j}x^{i}\beta
_{i,j}+\sum_{i=0}^{k-1}(-1)^{i+k-1}x^{i}\beta _{i,k-1}}{1-x-x(1-x)%
\sum_{j=1}^{k-2}\sum_{i=0}^{j}(-1)^{i+j}i\beta
_{i,j}-x\sum_{i=0}^{k-1}(-1)^{i+k-1}i\beta _{i,k-1}},
\end{equation*}%
where $\beta _{i,j}$ is defined in Theorem \ref{mth}.
\end{theorem}

For example, Theorem~\ref{mth22}, for $k=2,3,4,5,6$, gives the following
ordinary generating functions for the number of $(k,1)$-noncrossing
partitions of $[n]$:%
\begin{equation*}
\begin{array}{ll}
\sum\limits_{n\geq 0}\#\mathcal{N}_{2,1}(n)x^{n}=\frac{1-x}{1-2x}, &
\sum\limits_{n\geq 0}\#\mathcal{N}_{3,1}(n)x^{n}=\frac{1-3x+x^{2}}{%
(1-x)(1-3x)}, \\[14pt]
\sum\limits_{n\geq 0}\#\mathcal{N}_{4,1}(n)x^{n}=\frac{1-6x+9x^{2}-3x^{3}}{%
(1-x)(1-2x)(1-4x)}, & \sum\limits_{n\geq 0}\#\mathcal{N}_{5,1}(n)x^{n}=\frac{%
1-10x+32x^{2}-37x^{3}+11x^{4}}{(1-x)(1-2x)(1-3x)(1-5x)},%
\end{array}%
\end{equation*}%
\begin{equation*}
\begin{array}{l}
\sum\limits_{n\geq 0}\#\mathcal{N}_{6,1}(n)x^{n}=\frac{%
1-15x+81x^{2}-192x^{3}+189x^{4}-53x^{5}}{(1-x)(1-2x)(1-3x)(1-4x)(1-6x)}.%
\end{array}%
\end{equation*}%
The numbers of $(k,1)$-noncrossing partitions of $[n]$ with $k=2,3,4,5,6$
are given by
\begin{equation*}
\begin{array}{l}
\#\mathcal{N}_{2,1}(n)=2^{n-1}, \\[3pt]
\#\mathcal{N}_{3,1}(n)=\frac{1}{6}(3^{n}+3), \\[3pt]
\#\mathcal{N}_{4,1}(n)=\frac{1}{24}(4^{n}+6\cdot 2^{n}+8), \\[3pt]
\#\mathcal{N}_{5,1}(n)=\frac{1}{120}(5^{n}+10\cdot 3^{n}+20\cdot 2^{n}+45),
\\[3pt]
\#\mathcal{N}_{6,1}(n)=\frac{1}{720}(6^{n}+15\cdot 4^{n}+40\cdot
3^{n}+135\cdot 2^{n}+264).%
\end{array}%
\end{equation*}

\section{Proofs}

Let us denote by $F_{k}(x)$ the generating function for the number of
partitions in $\mathcal{P}_{k}(n)$:
\begin{equation*}
F_{k}(x)=\sum_{n\geq 0}\#\mathcal{P}_{k}(n)x^{n}.
\end{equation*}
Here, instead of dealing with recurrence relations with two
indices $n$ and $\ell$, as we mentioned in the introduction, we
deal with recurrence relations in terms of ordinary generating
functions with a single index $\ell$. Let us denote by $F_{k,\ell
}(x)$ the generating function for the number of partitions in
$\mathcal{P}_{k,\ell}(n)$:
\begin{equation*}
F_{k,\ell }(x)=\sum_{n\geq 0}\#\mathcal{P}_{k,\ell }(n)x^{n}
\end{equation*}
Here, for the case $\ell =0$ we have $F_{k,\ell }(x)=1$. Clearly,
$F_{k}(x)=\sum_{i\geq 0}F_{k,i}(x)$. Our main result is based on
the construction of linear recurrence relations with a single
index for the ordinary generating function $F_{k,\ell }(x)$. As we
will see later, since the recurrences contain the expression
$\sum_{i\geq\ell}F_{k,i}(x)$, we define for clarity
\begin{equation*}
G_{k,\ell }(x)=\sum_{i\geq \ell }F_{k,i}(x).
\end{equation*}
The expression $G_{k,\ell }(x)$ is the ordinary generating function
for the number of $j$-increasing partitions of $[n]$, with $j\geq
\ell$. It follows directly from the definitions that the generating
function $G_{k,\ell }(x)$ is well-defined, since
\begin{equation*}
G_{k,\ell }(x)=F_{k}(x)-\sum_{i=0}^{\ell -1}F_{k,i}(x).
\end{equation*}

In our first lemma, we find the recurrence relation for $F_{k,\ell }(x)$,
where $1\leq \ell \leq k-1$.

\begin{lemma}
\label{lem1} For all $1\leq \ell\leq k-1$,
\begin{equation*}
F_{k,\ell}(x)=\ell x G_{k,\ell}(x)+x^\ell.
\end{equation*}
\end{lemma}
\begin{proof}
Let $\pi =\pi _{1}\pi _{2}\cdots \pi _{n}\in \mathcal{P}_{k,\ell
}(n)$. If $n>\ell $ then $\pi _{\ell +1}\leq \ell $. Thus $\pi \in
\mathcal{P}_{k}(n)$ if and only if the reduced form of $\pi ^{\prime
}=\pi _{1}\cdots \pi _{\ell}\pi _{\ell +2}\cdots \pi _{n}$ is a
partition in $\cup _{j\geq \ell}\mathcal{P}_{k,j}(n-1)$. Therefore,
\begin{equation*}
F_{k,\ell }(x)=\ell x(F_{k,\ell }(x)+F_{k,\ell +1}(x)+\cdots )+x^{\ell
}=\ell xG_{k,\ell }(x)+x^{\ell },
\end{equation*}
where $x^{\ell }$ counts the unique $\ell$-increasing partition of
$[\ell]$ , namely $12\cdots\ell$, as required by the statement.
\end{proof}

The above observation together with the the definition of
$G_{k,\ell}(x)$ gives the following system:
\begin{equation}
\left\{
\begin{array}{ll}
F_{k,0}(x) & =1 \\
xF_{k,0}(x)+F_{k,1}(x) & =xF_{k}(x)+x \\
2xF_{k,0}(x)+2xF_{k,1}(x)+F_{k,2}(x) & =2xF_{k}(x)+x^{2} \\
\vdots  &  \\
(k-1)xF_{k,0}(x)+\cdots +(k-1)xF_{k,k-2}(x)+F_{k,k-1}(x) &
=(k-1)xF_{k}(x)+x^{k-1}.
\end{array}
\right. \label{eqa1}
\end{equation}
Next, we find an explicit formula for $F_{k,\ell}(x)$ in terms of
$F_{k}(x)$.

\begin{lemma}
\label{lem2} For all $1\leq \ell\leq k-1$,
\begin{equation*}
F_{k,\ell}(x)=\sum_{i=0}^\ell(-1)^{i+\ell}(ixF_k(x)+x^i)\beta_{i,\ell},
\end{equation*}
where $\beta_{\ell,\ell}=1$ and $\beta_{i,\ell}=\ell
x\prod_{j=i+1}^{\ell-1}(jx-1)$ for $i=0,1,\ldots,\ell-1$.
\end{lemma}

\begin{proof}
With the use of Cramer's Rule on \eqref{eqa1}, we obtain
\begin{equation*}
F_{k,\ell }(x)=\sum_{i=0}^{\ell }(-1)^{i+\ell }(ixF_{k}(x)+x^{i})\left\vert
\begin{array}{cccccc}
(i+1)x & 1 & 0 & \cdots  & 0 & 0 \\
(i+2)x & (i+2)x & 1 & \cdots  & 0 & 0 \\
\vdots  &  &  &  &  &  \\
(\ell -1)x & (\ell -1)x & (\ell -1)x & \cdots  & (\ell -1)x & 1 \\
\ell x & \ell x & \ell x & \cdots  & \ell x & \ell x%
\end{array}%
\right\vert .
\end{equation*}%
By making use of the formula
\begin{equation*}
\left\vert
\begin{array}{cccccc}
ax & 1 & 0 & \cdots  & 0 & 0 \\
(a+1)x & (a+1)x & 1 & \cdots  & 0 & 0 \\
\vdots  &  &  &  &  &  \\
(b-1)x & (b-1)x & (b-1)x & \cdots  & (b-1)x & 1 \\
bx & bx & bx & \cdots  & bx & bx
\end{array}
\right\vert =bx\prod_{j=a}^{b-1}(jx-1),
\end{equation*}
which holds by induction on $b\geq a$, we obtain
\begin{equation*}
F_{k,\ell }(x)=\sum_{i=0}^{\ell }(-1)^{i+\ell }(ixF_{k}(x)+x^{i})\beta
_{i,\ell },
\end{equation*}
as claimed.
\end{proof}

Now, before completing the proof of our main result, Theorem~\ref{mth}, let
us present two applications of Lemma~\ref{lem2}. The first one is the
ordinary generating function for the number of partitions in $\mathcal{Q}%
_{\ell }(n)$; the second one is the ordinary generating function for the
number of $(k,1)$-noncrossing partitions of $[n]$.

\subsection{Enumerating partitions in $\mathcal{Q}_\ell(n)$}
The formula of the ordinary generating function $I_{\ell }(x)$ for
the number of partitions in $\mathcal{Q}_{\ell }(n)$ can be obtained
as follows. From the definition of the set $Q_\ell(n)$ and from the
proof of Lemma \ref{lem2} for $\ell<k$, we obtain that the ordinary
generating function for the number of $m$-increasing partitions in
$Q_\ell(n)$ is given by
\begin{equation*}
I_{\ell,m}(x)=\sum_{i=0}^m(-1)^{i+m}(ixI_\ell(x)+x^i)\beta_{i,m}.
\end{equation*}
On the other hand, $I_\ell(x)=\sum_{m=0}^{\ell}I_{\ell,m}(x)$.
Combining these two equations, we obtain
\begin{equation*}
I_{\ell }(x)=\sum_{j=0}^{\ell }I_{\ell,j}(x)=1+\sum_{j=1}^{\ell
}\sum_{i=0}^{j}(-1)^{i+j}(ixI_{\ell }(x)+x^{i})\beta _{i,j}.
\end{equation*}
The solution of this equation gives a formula for $I_{\ell }(x)$, as stated
in Corollary \ref{se}.

\subsection{Enumerating $(k,1)$-noncrossing partitions of $[n]$}

Let $J_{k}(x)$ be the ordinary generating function for the number of $(k,1)$%
-noncrossing partitions of $[n]$, that is,%
\begin{equation*}
J_{k}(x)=\sum_{n\geq 0}\#\mathcal{N}_{k,1}(n)x^{n}.
\end{equation*}%
More generally, let $J_{k,\ell }(x)$ be the ordinary generating
function for the number of $(k,1)$-noncrossing $\ell $-increasing
partitions of $[n]$. Then, a similar argument as in the proof of
Lemma~\ref{lem2} gives that
\begin{equation}
J_{k,\ell }(x)=\sum_{i=0}^{\ell }(-1)^{i+\ell }(ixJ_{k}(x)+x^{i})\beta
_{i,\ell },\quad \ell =1,2,\ldots ,k-1,  \label{eqrr1}
\end{equation}%
with $J_{k,0}(x)=1$. On the other hand,
\begin{equation}
J_{k,\ell }(x)=x^{\ell +1-k}J_{k,k-1}(x),\quad \ell =k,k+1,k+2,\ldots .
\label{eqrr2}
\end{equation}%
To prove this observation, let $\pi $ be any $(k,1)$-noncrossing $\ell $%
-increasing partition of $[n]$ with $\ell \geq k$. Then $\pi _{1}\pi
_{2}\cdots \pi _{\ell }=12\cdots \ell $ and $\pi _{\ell +1}<\ell +1$. Since $%
\pi $ avoids $12\cdots k1$, then $\pi _{i}\notin \{1,2,\ldots ,\ell +1-k\}$
for all $i\geq \ell +1$. Thus, the number of $(k,1)$-noncrossing $\ell $%
-increasing partition of $[n]$ with $\ell \geq k$ is the same as the number
of $(k,1)$-noncrossing $(k-1)$-increasing partition of $[n-\ell -1+k]$. This
is equivalent to $J_{k,\ell }(x)=x^{\ell +1-k}J_{k,k-1}(x)$, for all $\ell
\geq k$. Therefore, using the fact that $J_{k}(x)=\sum_{\ell \geq
0}J_{k,\ell }(x)$, and the two equations \eqref{eqrr1} and \eqref{eqrr2}, we
can write
\begin{equation*}
J_{k}(x)=1+\sum_{\ell =0}^{k-1}\sum_{i=0}^{\ell }(-1)^{i+\ell
}(ixJ_{k}(x)+x^{i})\beta _{i,\ell }+\sum_{\ell \geq k}x^{\ell
+1-k}J_{k,k-1}(x).
\end{equation*}%
Again, by \eqref{eqrr1}, we have
\begin{equation*}
J_{k}(x)=1+\sum_{\ell =0}^{k-1}\sum_{i=0}^{\ell }(-1)^{i+\ell
}(ixJ_{k}(x)+x^{i})\beta _{i,\ell }+\frac{x}{1-x}%
\sum_{i=0}^{k-1}(-1)^{i+k-1}(ixJ_{k}(x)+x^{i})\beta _{i,k-1}.
\end{equation*}%
The solution of this equation gives a formula for $J_{k}(x)=\sum_{n\geq 0}\#%
\mathcal{N}_{k,1}(n)x^{n}$ as stated in Theorem \ref{mth22}.

\subsection{Proof of Theorem~\protect\ref{mth}}

We need some extra notation before completing the proof of Theorem~\ref{mth}%
. Let
\begin{equation*}
H_{k}(x,y)=\sum_{\ell =0}^{k-2}F_{k,\ell }(x)y^{\ell }
\end{equation*}%
and%
\begin{equation*}
F_{k}(x,y)=\sum_{\ell \geq 0}F_{k,\ell }(x)y^{\ell }.
\end{equation*}%
From Lemma~\ref{lem2}, we can observe that
\begin{equation*}
H_{k}(x,y)=\sum_{\ell =0}^{k-2}y^{\ell }\left( \sum_{i=0}^{\ell
}(-1)^{i+\ell }(ixF_{k}(x,1)+x^{i})\beta _{i,\ell }\right)
\end{equation*}%
and that
\begin{equation*}
F_{k}(x,1)=F_{k}(x),
\end{equation*}%
where $\beta _{\ell ,\ell }=1$ and $\beta _{i,\ell }=\ell
x\prod_{j=i+1}^{\ell -1}(jx-1)$, for $i=0,1,\ldots ,\ell -1$. Now, let us
focus on the generating functions $F_{k,\ell }(x)$, where $\ell \geq k-1$.

\begin{lemma}
\label{lem3} For all $j\geq0$,
\begin{equation*}
F_{k,k-1+j}(x)=x^{k-1+j}+%
\sum_{i=0}^{j-1}x^{j+1-i}G_{k,k-1+i}(x)+(k-1)xG_{k,k-1+j}(x).
\end{equation*}
\end{lemma}

\begin{proof}
The case $j=0$ holds on the basis of Lemma~\ref{lem1}. Let us assume that $%
j\geq 1$. Let $\pi =\pi _{1}\pi _{2}\cdots \pi _{n}$ be any partition in $%
\mathcal{P}_{k,k-1+j}(n)$ such that $\pi _{k+j}=i\leq k-1+j$. We want the
equation of the generating function for the number of partitions in $\mathcal{P}%
_{k,k-1+j}(n)$, namely $F_{k,k-1+j}(x)$. Let us consider the following two
cases:

\begin{itemize}
\item If $1\leq i\leq j$ then $\pi $ is such that $\pi _{p}\notin
\{i+1,i+2,\ldots ,j+1\}$, where $p\geq k+1+j$. Thus, the contribution of
this case is
\begin{equation*}
x^{j+1-i}(F_{k,k-1+i}(x)+F_{k,k+i}(x)+\cdots )=x^{j+1-i}G_{k,k-1+i}(x).
\end{equation*}

\item If $j+1\leq i\leq k-1+j$ then $\pi $ satisfies the above conditions if
and only if the reduced form of $\pi _{1}\cdots \pi _{k-1+j}\pi
_{k+1+j}\cdots \pi _{n}$ is a partition in $\cup _{i\geq 0}\mathcal{P}%
_{k,k-1+j+i}(n-1)$. Thus, the contribution of this case is
\begin{equation*}
x(F_{k,k-1+j}(x)+F_{k,k+j}(x)+\cdots )=xG_{k,k-1+j}(x).
\end{equation*}
\end{itemize}

Putting together the above cases, $i=1,2,\ldots ,k-1+j$, we obtain that
\begin{equation*}
F_{k,k-1+j}(x)=x^{k-1+j}+%
\sum_{i=0}^{j-1}x^{j+1-i}G_{k,k-1+i}(x)+(k-1)xG_{k,k-1+j}(x),
\end{equation*}%
where $x^{k-1+j}$ counts the unique $(k-1+j)$-increasing partitions of $%
[k-1+j]$, namely $12\cdots (k-1+j)$.
\end{proof}

Now we have a formula for the generating function $F_{k}(x,y)$:

\begin{proposition}
\label{pro1} We have
\begin{equation*}
\begin{array}{l}
\left(1+\frac{x^2y^2}{(1-y)(1-xy)}+\frac{(k-1)xy}{1-y}
\right)(F_k(x,y)-H_k(x,y)) \\[4pt]
\qquad\qquad\qquad\qquad\qquad\qquad=\frac{(xy)^{k-1}}{1-xy}+\frac{y^{k-1}}{
1-y}\left(\frac{x^2y}{1-xy}+(k-1)x\right)(F_k(x,1)-H_k(x,1)).
\end{array}
\end{equation*}
\end{proposition}
\begin{proof}
Lemma~\ref{lem1} together with Lemma~\ref{lem3} give
\begin{equation*}
\begin{array}{l}
F_{k}(x,y) \\
=\frac{1}{1-xy}+\sum\limits_{j=1}^{k-2}jxG_{k,j}(x)y^{j}+(k-1)x\sum%
\limits_{j\geq k-1}G_{k,j}(x)y^{j}+\frac{x^{2}y}{1-xy}\sum\limits_{j\geq
k-1}G_{k,j}(x)y^{j} \\[11pt]
=\frac{1}{1-xy}+\sum\limits_{j=1}^{k-2}(F_{k,j}(x)-x^{j})y^{j} \\
\ \ \ +x\left( \frac{xy}{1-xy}+k-1\right) \left( \frac{y^{k-1}}{1-y}%
(F_{k}(x,1)-H_{k}(x,1))-\frac{y}{1-y}(F_{k}(x,y)-H_{k}(x,y))\right) ,%
\end{array}
\end{equation*}
which is equivalent to
\begin{equation*}
\begin{array}{l}
\left( 1+\frac{x^{2}y^{2}}{(1-y)(1-xy)}+\frac{(k-1)xy}{1-y}\right)
(F_{k}(x,y)-H_{k}(x,y)) \\[4pt]
\qquad \qquad \qquad \qquad \qquad \qquad =\frac{(xy)^{k-1}}{1-xy}+\frac{%
y^{k-1}}{1-y}\left( \frac{x^{2}y}{1-xy}+(k-1)x\right)
(F_{k}(x,1)-H_{k}(x,1)),
\end{array}
\end{equation*}
as claimed.
\end{proof}

The functional equation in the statement of Proposition~\ref{pro1} can be
solved systematically using the \emph{kernel method} \emph{technique }(see
\cite{B02}). Let
\begin{equation*}
y=y_{k}=\frac{1-(k-2)x-\sqrt{(1-kx)^{2}-4x^{2}}}{2x(1-(k-2)x)}
\end{equation*}
be one of the roots of the equation
$1+\frac{x^{2}y^{2}}{(1-y)(1-xy)}+\frac{(k-1)xy}{1-y}$. Then
Proposition~\ref{pro1} gives
\begin{equation*}
F_{k}(x,1)-H_{k}(x,1)=\frac{x^{k-1}y_{k}}{1-xy_{k}}.
\end{equation*}
Therefore, Lemma~\ref{lem2} gives
\begin{equation*}
F_{k}(x,1)-\sum_{j=0}^{k-2}\sum_{i=0}^{j}(-1)^{i+j}(ixF_{k}(x,1)+x^{i})\beta
_{i,j}=\frac{x^{k-1}y_{k}}{1-xy_{k}},
\end{equation*}
which implies our main result, Theorem~\ref{mth}.

\medskip

\textbf{Acknowledgment}. Part of this work has been carried out while the
second author was visiting the Center for Combinatorics of Nankai
University. The financial support of the Center is gratefully acknowledged.
Also, the authors would like to thank the anonymous referees for a number of
valuable comments, which helped to improve the presentation of this paper.



\begin{thebibliography}{99}
\bibitem{B02} \textsc{C.~Banderier, M.~Bousquet-M{\'{e}}lou, A.~Denise,
P.~Flajolet, D.~Gardy, and D.~Gouyou-Beauchamps}, \newblock Generating
functions for generating trees, \newblock Formal Power Series and Algebraic
Combinatorics (Barcelona, 1999), \emph{Discr. Math.} \textbf{246}:1-3 (2002)
29--55.

\bibitem{C1} \textsc{William Y. C. Chen, Eva Y.P. Deng, and Rosena R.X. Du},
Reduction of $m$-regular noncrossing Partitions, \emph{Europ. J. Combin.}
\textbf{26}:2 (2005) 237--243.

\bibitem{C2} \textsc{William Y. C. Chen, Eva Y.P. Deng, Rosena R.X. Du, R.
P. Stanley, and Catherine H. Yan}, Crossings and Nestings of Matchings and
Partitions, \emph{Trans. Amer. Math. Soc.}, to appear. arXiv:
math.CO/0501230.


\bibitem{K1} \textsc{M. Klazar}, On abab-free and abba-free set partitions,
\emph{Europ. J. Combin.} \textbf{17} (1996) 53--68.

\bibitem{K2} \textsc{M. Klazar}, On trees and noncrossing partitions, \emph{%
Discr. Appl. Math.} \textbf{82} (1998) 263--269.


\bibitem{JS} \textsc{N.J.A. Sloane}, The Online Encyclopedia of Integer
Sequences, \emph{www.research.att.com/$\sim$njas/ sequences/}.

\bibitem{SA} \textsc{B. E. Sagan}, Pattern avoidance in set partitions,
preprint.

\bibitem{S1} \textsc{R. P. Stanley}, Enumerative Combinatorics, Vol. 1,
Cambridge University Press, Cambridge, UK, 1996.

\bibitem{WW} \textsc{M. Wachs and D. White}, p,q-Stirling numbers and set
partition statistics, \emph{J. Combin. Theory, Series A}, \textbf{56}:1
(1991) 27--46.
\end{thebibliography}
\end{document}